\documentclass[a4paper]{article}
\usepackage[margin=1.2in]{geometry}

\usepackage{amssymb}
\usepackage{amsmath,amsthm}
\usepackage{graphicx}
\usepackage{color}
\usepackage{booktabs,tabularx, colortbl} 
\usepackage{authblk}
\usepackage{setspace}

\definecolor{dkgreen}{rgb}{0,0.6,0} 
\definecolor{gray}{rgb}{0.5,0.5,0.5}
\definecolor{mauve}{rgb}{0.58,0,0.82}

\newcommand{\bbR}{\mathbb{R}}

\newcommand{\bbZ}{\mathbb{Z}}
\newcommand{\mO}{\mathcal{O}}

\newtheorem{thm}{Theorem}
\newtheorem{lem}{Lemma}

\newtheorem{prop}{Proposition}

\newtheorem{ass}{Assumption}
\newtheorem{obs}{Observation}
\newtheorem{rem}{Remark}
\newtheorem{clm}{Claim}

\begin{document}
\title{Improving bounds on the diameter of a polyhedron in high dimensions}
\author[]{Noriyoshi Sukegawa\thanks{sukegawa@ise.chuo-u.ac.jp}}
\affil[]{{\normalsize Department of Information and System Engineering, 
Faculty of Science and Engineering, Chuo University, 
1-13-27 Kasuga, Bunkyo-ku, Tokyo 112-8551, Japan}}

\maketitle

\begin{abstract}
In 1992, Kalai and Kleitman proved that the diameter of a $d$-dimensional polyhedron with $n$ facets is at most $n^{2+\log_2 d}$. 
In 2014, Todd improved the Kalai-Kleitman bound to $(n-d)^{\log_2 d}$. 
We improve the Todd bound to $(n-d)^{-1+\log_2 d}$ for $n \ge d \ge 7$, $(n-d)^{-2+\log_2 d}$ for $n \ge d \ge 37$, and $(n-d)^{-3+\log_2 d+O\left(1/d\right)}$ for $n \ge d \ge 1$. 
\\
\vspace{.05in}
\\
{\bf keyword:}~
Diameter; polyhedra; high dimension; computer-assisted method
\end{abstract}

\section{Introduction}\label{sec:intro}

The \emph{diameter} $\delta(P)$ of a polyhedron $P$ is the smallest integer $k$ such that every pair of vertices of $P$ can be connected by a path using at most $k$ edges of $P$. 
The diameter is a fundamental feature of a polyhedron and is closely related to the theoretical complexity of the simplex algorithm; 
the number of pivots needed, in the worst case, by the simplex algorithm to solve a linear program on a polyhedron $P$ is bounded from below by $\delta(P)$. 

One of the outstanding open problems in the areas of polyhedral combinatorics and operations research is to understand 
the behavior of $\Delta(d,n)$, the maximum possible diameter of a $d$-dimensional polyhedron with $n$ facets. 
In 1957, Warren M. Hirsch asked whether $\Delta(d,n) \le n-d$. 
While this inequality was shown to hold for $d \le 3$~\cite{Kl64,Kl65,Kl66}, Klee and Walkup~\cite{KlWa67} disproved it for unbounded polyhedra when $d \ge 4$ in 1967, 
and Santos~\cite{Sa12} finally  disproved it for bounded polyhedra, i.e., for polytopes, in 2012. 
Santos' lower bound, later refined by Matschke, Santos, and Weibel \cite{MaSaWe12}, however, violates $n-d$ by only $5$ percent. 
For the history of the Hirsch conjecture, see \cite{Sa13}. 

The first subexponential upper bound on $\Delta(d,n)$ is due to Kalai and Kleitman~\cite{KaKl92} 
who proved in 1992 that $\Delta(d,n)$ is at most $n^{2+\log_2 d}$. 
The key ingredient for their proof is a recursive inequality on $\Delta(d,n)$, which we call the {\em Kalai-Kleitman inequality}. 
The Kalai-Kleitman inequality was later extended to more general settings such as connected layer families by Eisenbrand et al.~\cite{Ei10}, and subset partition graphs by Gallagher and Kim~\cite{GaKi14}. 
For the corresponding lower bounds, we refer to \cite{Ei10,Ki14}. 

Refining Kalai and Kleitman's approach, in \cite{To14}, Todd showed in 2014 that $\Delta(d,n)\le (n-d)^{\log_2 d}$ for $n \ge d \ge 1$. 
The Todd bound is tight for $d \le 2$ and coincides with the true value $\Delta(d,d)$, i.e., $0$, when $n=d$. 
Sukegawa and Kitahara~\cite{SuKi15} slightly improved the Todd bound to $(n-d)^{\log_2(d-1)}$ for $n \ge d \ge 3$. 
We note that their bound is no longer valid for $d \le 2$, however, it coincides with the Hirsch bound of $n-d$, and is tight for $d=3$. 
On the other hand, Gallagher and Kim~\cite{GaKi16} proved that the same bound holds for the diameter of normal simplicial complexes, and also improved it for polytopes.  

\subsection{Main results} 

In this paper, we improve the Todd bound in high dimensions as follows: 
\begin{thm}\label{thm:main}
$\:$
\begin{description}\setlength{\itemsep}{0pt}
\item[$(a)$]$\Delta(d,n)\le (n-d)^{\log_2\left( d/2 \right)} = (n-d)^{-1+\log_2 d}$ for $n \ge d \ge 7$,
\item[$(b)$]$\Delta(d,n)\le (n-d)^{\log_2\left( d/4 \right)} = (n-d)^{-2+\log_2 d}$ for $n \ge d \ge 37$, and 
\item[$(c)$]$\Delta(d,n)\le (n-d)^{\log_2\left( 16 + d/8 \right)} = (n-d)^{-3+\log_2 d+\mO\left(1/d\right)}$ for $n \ge d \ge 1$. 
\end{description}
\end{thm}

\noindent
Inequalities $(a)$ and $(b)$ hold for, respectively, $d \ge 7$ and $d \ge 37$, and improve the Todd bound by, respectively, one and two orders of magnitude. 
Inequality $(c)$ holds for any $d$, and improves the Todd bound for $d \ge 19$.  
Note that $\log_2\left(16+\frac{d}{8}\right) = \log_2(d)-3+O\left(\frac{1}{d}\right)$ since $\log_e (1+x) \le x$ for $x \ge 0$. 
Thus, Inequality $(c)$ improves the Todd bound by roughly three orders of magnitude for sufficiently large $d$. 

\subsection{Our approach} 

As in~\cite{KaKl92,SuKi15,To14}, each inequality stated in Theorem~\ref{thm:main} will be proved via an induction on $d$ based on the Kalai-Kleitman inequality. 
In contrast to \cite{KaKl92,SuKi15,To14}, we introduce a way of strengthening Todd's analysis for the inductive step in high dimensions. 
In this approach, on the other hand, we need to check a large number of pairs $(d,n)$ for the base case. 
To address this issue, we devise a computer-assisted method which is based on two previously known upper bounds on $\Delta(d,n)$: 
\begin{description}\setlength{\itemsep}{0pt}
\item[$(i)$] $\tilde{\Delta}(d,n)$, an {\em implicit} upper bound on $\Delta(d,n)$ computed recursively from the Kalai-Kleitman inequality, 
\item[$(ii)$] the generalized Larman bound implying $\Delta(d,n) \le 2^{d-3}n$. 
\end{description}
\noindent
The Larman bound of $2^{d-3}n$ was originally proved for bounded polyhedra~\cite{La70}, and improved to $\frac{2n}{3}2^{d-3}$ by Barnette~\cite{Ba74}. 
Considering a more generalized setting, Eisenbrand et al.~\cite{Ei10} proved a bound of $2^{d-1}n$ in 2010, 
before Labb{\'e}, Manneville, and Santos~\cite{LaMaSa15} established in 2015 an upper bound on the diameter of simplicial complexes implying $\Delta(d,n) \le 2^{d-3}n$. 

\subsection{Related work} 

It should be noted that although this paper deals with only the two parameters $d$ and $n$, i.e., the dimension and the number of facets of a polyhedron, there have been studies on other parameters. 

A well-known example is the maximum integer coordinate of \emph{lattice polytopes}. 
In \cite{KlOn92}, Kleinschmidt and Onn proved that the diameter of a lattice polytope whose vertices are drawn from $\{0,1,\ldots,k\}^d$ is at most $kd$.
This is an extension of Naddef~\cite{Na89} showing that the diameter of a $0$-$1$ polytope is at most $d$. 
In 2015, Del Pia and Michini~\cite{DeMi15} improved the Kleinschmidt-Onn bound to $kd-\lceil \frac{d}{2} \rceil$ for $k \ge 2$ and showed that it is tight for $k=2$, 
before Deza and Pournin~\cite{DePo16} further improved the bound to $kd-\lceil \frac{2d}{3} \rceil - (k-3)$ for $k \ge 3$. 
On the other hand, considering Minkowski sums of primitive lattice vectors, in \cite{DeMaOn15}, Deza, Manoussakis, and Onn provided a lower bound of $\lfloor \frac{(k+1)d}{2} \rfloor$ for $k < d$. 

Another well-studied parameter would be $\Delta_A$ which is defined as the largest absolute value of a subdeterminant of the constraint matrix $A$ associated to a polyhedron. 
Bonifas et al.~\cite{Bo14} strengthened and extended the Dyer and Frieze upper bound~\cite{DyFr94} holding for totally unimodular case; i.e., when $\Delta_A=1$. 
Complexity analyses based on $\Delta_A$ for the shadow vertex algorithm and the primal-simplex based Tardos' algorithm were proposed by 
Dadush and H\"{a}hnle~\cite{DaHa14}, and Mizuno, Sukegawa, and Deza~\cite{MiSuDe14,MiSuDe15}, respectively.

We also note that there are studies that attempt to understand the behavior of $\Delta(d,n)$ when the number of facets is sufficiently large. 
Gallagher and Kim~\cite{GaKi16} provided an upper bound on the diameter of a normal simplicial complex and showed the {\em tail-polynomiality}; 
more specifically, they showed that the diameter is bounded from above by a polynomial in $n$ when $n$ is sufficiently large. 
An alternative simpler proof for such tail-polynomial upper bounds can be found in Mizuno and Sukegawa~\cite{MiSu16}. 
In contrast, in this paper, we assume that $d$ is large, and try to utilize this assumption to strengthen the previous results. 

\section{Preliminaries} 

A \emph{polyhedron} $P \subseteq \bbR^d$ is an intersection of a finite number of closed halfspaces, and $\dim(P)$ denotes the dimension of the affine hull of $P$. 
For a polyhedron $P$, an inequality $a^\top x \le \beta$ is said to be \emph{valid} for $P$ if it is satisfied by every $x \in P$. 
We say that $F$ is a \emph{face} of $P$ if there is a valid inequality $a^\top x \le \beta$ for $P$ which satisfies $F=P \cap \{ x \in \bbR^d: a^\top x = \beta \}$. 
In particular, $0$-, $1$-, and $(\dim(P)-1)$-dimensional faces are, respectively, referred to as \emph{vertices}, \emph{edges}, and \emph{facets}.   

The \emph{diameter} $\delta(P)$ of a polyhedron $P$ is the smallest integer $k$ such that every pair of vertices of $P$ can be connected by a path using at most $k$ edges of $P$.
In this paper, we are concerned with upper bounds on $\Delta(d,n)$, the maximum possible diameter of a $d$-dimensional polyhedron with $n$ facets. 
Lemma~\ref{thm:KaKl-ineq} states the {\em Kalai-Kleitman inequality} on which our approach is based.
\begin{lem}[Kalai-Kleitman inequality~\cite{KaKl92}]\label{thm:KaKl-ineq}
For $\left\lfloor \frac{n}{2} \right\rfloor \ge d \ge 2$, 
\begin{align*}\label{eqn:KaKl}
\Delta(d, n) \le \Delta(d-1, n-1) + 2 \Delta\left(d, \left\lfloor \frac{n}{2} \right\rfloor \right) + 2.
\end{align*}\noindent
\end{lem}

\subsection{Basic idea of our proof} 

We consider upper bounds of the form: \[f_{\alpha,\beta}(d,n)=(n-d)^{\log_2(\beta+d/\alpha)}, \]
where $(\alpha, \beta) \in S = \{ (\alpha,\beta) \in \bbZ^2: \alpha > 0, \ \beta \geq 0 \}$ is a pair of integers controlling the {\em quality} of upper bounds.  
Note that the Todd bound is $f_{1,0}(d,n)$. 
The upper bounds appearing in Inequalities $(a)$, $(b)$, and $(c)$ stated in Theorem~\ref{thm:main} correspond, respectively, to $f_{2,0}(d,n)$, $f_{4,0}(d,n)$, and $f_{8,16}(d,n)$. 

As mentioned earlier, we prove Inequalities $(a)$, $(b)$, and $(c)$ stated in Theorem~\ref{thm:main} via an induction on $d$ based on the Kalai-Kleitman inequality. 
The following lemma is the key ingredient for the inductive step; see Section~\ref{sec:correctness} for a proof. 
\begin{lem}\label{lem:validInductive}
If $(\alpha, \beta) \in S$, 
then there exists $d(\alpha,\beta)$  such that $d \ge d(\alpha,\beta)$, $n \ge 2d$, and $n \ge d+2^{2\alpha+1}$ imply
\begin{align}\label{eqn:inductive}
f_{\alpha,\beta}(d-1,n-1)+2f_{\alpha,\beta}\left(d, \left\lfloor \frac{n}{2} \right\rfloor \right)+2 \le f_{\alpha,\beta}(d,n). 
\end{align}\noindent
\end{lem}

\subsubsection{Inductive step}\label{induct} 

Assume $d \ge d(\alpha,\beta)$ and \texttt{P$_{d-1}$}: $\Delta(d-1,n) \le f_{\alpha,\beta}(d-1,n)$ for $n \ge d-1$, as the induction hypothesis on $d$. 
In what follows, by induction on $n$, we prove \texttt{P$_{d}$}: $\Delta(d,n) \le f_{\alpha,\beta}(d,n)$ for $n \ge d$. 
First, let us consider the case $n<2d$. 
In this case, the claim, i.e., the desired inequality $\Delta(d,n) \le f_{\alpha,\beta}(d,n)$, follows from the following fundamental proposition; for a proof, see, e.g., \cite{To14}. 
\begin{prop}\label{prop:n<2d}
$\Delta(d,n) \le \Delta(d-1,n-1)$ for $n < 2d$. 
\end{prop}\noindent
From Proposition~\ref{prop:n<2d} and \texttt{P$_{d-1}$}, for $n < 2d$, 
\begin{align*}
\Delta(d,n) 	&\le \Delta(d-1,n-1) \le f_{\alpha,\beta}(d-1,n-1) \le f_{\alpha,\beta}(d,n), 
\end{align*}
where the last inequality follows since $\alpha>0$. 

Now, suppose that $n \ge 2d$. 
First, let us consider the case $n < d+2^{2\alpha+1}$. 
We observe that the number of integers $n$ satisfying the condition, i.e., $2d \le n < d+2^{2\alpha+1}$, is finite for fixed $d$, and becomes zero for $d \ge 2^{2\alpha+1}$. 
We therefore verify $\Delta(d,n) \le f_{\alpha,\beta}(d,n)$ for these pairs as a part of the base case. 
Next, let us consider the case $n \ge d+2^{2\alpha+1}$. 
In this case, we apply the Kalai-Kleitman inequality to yield
\begin{align*}\label{eq:mid}
\displaystyle
\Delta(d,n) 	&\le \Delta(d-1,n-1) + 2\Delta\left(d, \left\lfloor \frac{n}{2} \right\rfloor \right) +2\\
& \le f_{\alpha,\beta}(d-1,n-1)+2f_{\alpha,\beta}\left(d, \left\lfloor \frac{n}{2} \right\rfloor \right)+2,
\end{align*}
where the second inequality follows from the induction hypotheses on $d$ and $n$. 
Note that Lemma~\ref{lem:validInductive} applies to this case, which yields the desired inequality $\Delta(d,n) \le f_{\alpha,\beta}(d,n)$. 

\subsubsection{Base case} 

\begin{prop}\label{prop:0}
Let $(\alpha, \beta) \in S$. If there exists $l$ satisfying 
\begin{description}\setlength{\itemsep}{0pt}
\item[{\rm $(B_0)$}] case $d=l$: $\Delta(d,n) \le f_{\alpha,\beta}(d,n)$ for $n \ge d$, 
\item[{\rm $(B_1)$}] case $l < d < d(\alpha,\beta)$: $\Delta(d,n) \le f_{\alpha,\beta}(d,n)$ for $n \ge 2d$, 
\item[{\rm $(B_2)$}] case $d(\alpha,\beta) \le d < 2^{2\alpha+1}$: $\Delta(d,n) \le f_{\alpha,\beta}(d,n)$ for $n$ with $2d \le n < d+2^{2\alpha+1}$, 
\end{description}
then $\Delta(d,n) \le f_{\alpha,\beta}(d,n)$ for $d \ge l$
\begin{proof}
By similar arguments used in the inductive step in Section~\ref{induct}, 
$\Delta(d,n) \le f_{\alpha,\beta}(d,n)$ for $l < d < d(\alpha,\beta)$ if $(B_0)$ and  $(B_1)$ hold.
Similarly, $\Delta(d,n) \le f_{\alpha,\beta}(d,n)$ for $d \ge d(\alpha,\beta)$ if $(B_0)$, $(B_1)$, and $(B_2)$ hold.
\end{proof}
\end{prop}

\noindent
In this study, we devise a computer-assisted method to test whether $(B_0)$, $(B_1)$, and $(B_2)$ hold or not in a finite process. 
To this end, we 
\begin{description}\setlength{\itemsep}{0pt}
\item[{\rm ($I_1$)}] 
make the number of pairs $(d,n)$ to be checked in $(B_0)$ and $(B_1)$ finite, and 
\item[{\rm ($I_2$)}] 
establish an upper bound $\tilde{\Delta}(d,n)$ on $\Delta(d,n)$ which enables us to ensure $\Delta(d,n) \le f_{\alpha,\beta}(d,n)$ 
via the relationship $\Delta(d,n) \le \tilde{\Delta}(d,n) \le f_{\alpha,\beta}(d,n)$.  
\end{description}
\noindent
When $d \le 3$, one can set, for example, $\tilde{\Delta}(d,n) := n-d$ in {\rm ($I_2$)}. 
However, for large $d$, previously known upper bounds on $\Delta(d,n)$, including the Todd bound, are of course greater than $f_{\alpha,\beta}(d,n)$, 
and therefore cannot be used for deriving the desired inequality, i.e., $\tilde{\Delta}(d,n) \le f_{\alpha,\beta}(d,n)$. 
This is the reason why we need a computer-assisted method. 

\subsubsection{Strategy to {\rm ($I_1$)}}

We first explain our strategy to {\rm ($I_1$)}, i.e., how to make the number of pairs $(d,n)$ to be checked in $(B_0)$ and $(B_1)$ finite. 
\begin{ass}\label{ass:super-linear}
The choice of $(\alpha,\beta) \in S$ is such that 
$f_{\alpha,\beta}(d,n)=(n-d)^{\log_2(\beta+d/\alpha)}$ is superlinear in $n$ for fixed $d$ when $d \ge l$; 
i.e., $\alpha$ and $\beta$ satisfy $\log_2\left(\beta+\frac{d}{\alpha}\right)>1$ for $d \ge l$. 
\end{ass}

\begin{obs}\label{obs:Larman}
Suppose that $(\alpha,\beta) \in S$ satisfies Assumption~\ref{ass:super-linear}. 
For fixed $d$, if we let $n_L(d)$ be the smallest integer $n$ such that $2^{d-3}n \le f_{\alpha,\beta}(d,n)$ holds, 
then for $n \ge n_L(d)$, \[\Delta(d,n) \le 2^{d-3}n \le f_{\alpha,\beta}(d,n).\]   
\begin{proof}
Direct consequence of the generalized Larman bound. 
\end{proof}
\end{obs}
Thus, with Assumption~\ref{ass:super-linear}, if 
\begin{description}\setlength{\itemsep}{0pt}
\item[{\rm $(B'_0)$}] case $d=l$: $\Delta(d,n) \le f_{\alpha,\beta}(d,n)$ for $d \le n \le n_L(d)$, 
\item[{\rm $(B'_1)$}] case $l < d < d(\alpha,\beta)$: $\Delta(d,n) \le f_{\alpha,\beta}(d,n)$ for $2d \le n \le n_L(d)$, 
\end{description}
are satisfied, then $(B_0)$ and $(B_1)$ are satisfied. 
\begin{prop}\label{prop:0'}
Let $(\alpha, \beta) \in S$. 
If there exists $l$ satisfying Assumption~\ref{ass:super-linear}, and $(B'_0)$, $(B'_1)$, and $(B_2)$, 
then $\Delta(d,n) \le f_{\alpha,\beta}(d,n)$ for $d \ge l$. 
\end{prop}

\begin{table}[tb]
\caption{The pairs $(d,n)$ for which $\Delta(d,n) \le f_{\alpha,\beta}(d,n)$ should be ensured for Proposition~\ref{prop:0'}\label{fig:base}}
\centering
\begin{tabular}{lccccccccccccccccccccccccccccccccccccc}\toprule
&		&\multicolumn{12}{c}{$n-2d$}\\ \cmidrule(r){3-14}
Case		&$d$	&	0	&	1	&	2	&\multicolumn{8}{c}{$\ldots$}&		\\	\midrule
&$l+1$				&	-	&	$\triangleleft$	&		&		&		&		&		&		&		&		&		&		\\ 
&					&	-	&	-	&	$\triangleleft$	&		&		&		&		&		&		&		&		&		\\ 
$(B'_1)$&				&	-	&	-	&	-	&	-	&	$\triangleleft$	&		&		&		&		&		&		&		\\ 
&					&	-	&	-	&	-	&	-	&	-	&	-	&	-	&	$\triangleleft$	&		&		&		&		\\ 
&					&	-	&	-	&	-	&	-	&	-	&	-	&	-	&	-	&	-	&	-	&	-	&	$\triangleleft$	\\ \arrayrulecolor{gray} \midrule \arrayrulecolor{black} 
&$d(\alpha,\beta)$		&	-	&	-	&	-	&	-	&	-	&	-	&	$\circ$	&		&		&		&		&		\\ 
&					&	-	&	-	&	-	&	-	&	-	&	$\circ$	&		&		&		&		&		&		\\ 
&					&	-	&	-	&	-	&	-	&	$\circ$	&		&		&		&		&		&		&		\\ 
$(B_2)$&				&	-	&	-	&	-	&	$\circ$	&		&		&		&		&		&		&		&		\\ 
&					&	-	&	-	&	$\circ$	&		&		&		&		&		&		&		&		&		\\ 
&					&	-	&	$\circ$	&		&		&		&		&		&		&		&		&		&		\\  
&$2^{2\alpha+1}$		&	$\circ$	&		&		&		&		&		&		&		&		&		&		&		\\	\bottomrule \\
\end{tabular}
\end{table}

\noindent
The total number of pairs $(d,n)$ to be checked in $(B'_0)$, $(B'_1)$, and $(B_2)$ is {\em finite} as illustrated in Table~\ref{fig:base}. 
In the table, we assume that we have already found a dimension $l$ satisfying $(B'_0)$ and that $l+1<d(\alpha,\beta)$. 
If $l+1 \ge d(\alpha,\beta)$, then there is no pair $(d,n)$ to be checked for $(B'_1)$. 
Also, if $d(\alpha,\beta) \ge 2^{2\alpha+1}$, then the table will be much simpler since there is no pair $(d,n)$ to be checked for $(B'_1)$ and $(B_2)$. 
The pairs $(d,n)$ with $n<2d$ are omitted as the desired inequalities hold inductively. 
The meanings of the symbols are as follows:
\begin{itemize}\setlength{\itemsep}{0pt}
\item[-]: corresponds to a pair $(d,n)$ to which $\tilde{\Delta}(d,n) \le f_{\alpha,\beta}(d,n)$ must be ensured,
\item[$\triangleleft$]: corresponds to a pair $(d,n)$ with $n=n_L(d)$, and 
\item[$\circ$]: corresponds to a pair $(d,n)$ with $n=d+2^{2\alpha+1}$. 
\end{itemize}

\begin{rem}[How to compute $n_L(d)$]
In practice, we do not need to compute the value of $n_L(d)$ in advance. 
It suffices to check if $2^{d-3}n \le f_{\alpha,\beta}(d,n)$ for $n=2d,2d+1,\ldots$, for each fixed $d$. 
If $2^{d-3}n \le f_{\alpha,\beta}(d,n)$ holds for the first time for some pair $(d,n')$, then $n_L(d)=n'$. 
\end{rem}

\subsubsection{Strategy to {\rm ($I_2$)}}
We now explain our strategy to {\rm ($I_2$)}, i.e., how to establish an upper bound $\tilde{\Delta}(d,n)$ on $\Delta(d,n)$ which enables us to ensure $\Delta(d,n) \le f_{\alpha,\beta}(d,n)$ 
via the relationship $\Delta(d,n) \le \tilde{\Delta}(d,n) \le f_{\alpha,\beta}(d,n)$. 
We note that our strategy is based on Todd~\cite{To14}. 
Specifically, we define $\tilde{\Delta}(d,n)$ as a value recursively computed via:
\begin{align*}
\tilde{\Delta}(d, n)=
\begin{cases}
\ n-3 &\mbox{if }d = 3\mbox{ and }n \ge d,\\
\ \tilde{\Delta}(d-1, n-1)&\mbox{if }d>3\mbox{ and }d \le n<2d,\\
\ \tilde{\Delta}(d-1, n-1) + 2 \tilde{\Delta}\left(d, \left\lfloor \frac{n}{2} \right\rfloor \right) + 2 &\mbox{if }d>3\mbox{ and }n \ge 2d. \\
\end{cases}
\end{align*}
For example, 
\begin{align*}
\tilde{\Delta}(5, 13) &=\tilde{\Delta}(4, 12)+2\tilde{\Delta}(5, 6) +2\\
&=\left[ \tilde{\Delta}(3, 11)+ 2\tilde{\Delta}(4, 6) + 2 \right]+2\tilde{\Delta}(4, 5)+2\\
&=\left[ (11-3)+ 2\tilde{\Delta}(3, 5) + 2 \right]+2\tilde{\Delta}(3,4)+2\\
&=\left[ 10+ 2(5-3) \right]+2(4-3)+2=18. 
\end{align*}
\noindent
Then, by the validity of the Kalai-Kleitman inequality, Proposition~\ref{prop:n<2d}, and the correct inequality $\Delta(3,n) \le n-3$, 
we have $\Delta(d,n) \le \tilde{\Delta}(d, n)$ for every pair $(d,n)$ with $n \ge d \ge 3$. 
Therefore, for the pairs $(d,n)$ indicated by ``-'' in Table~\ref{fig:base}, we check if $\tilde{\Delta}(d, n) \le f_{\alpha,\beta}(d,n)$ instead of $\Delta(d, n)\le f_{\alpha,\beta}(d,n)$. 

\section{Proof Method}\label{sec:method} 

The section is devoted to the detailed description of the computer-assisted method for verifying the base case. 
We show our code in the programming language C and its execution results in Appendix~\ref{append:A}. 
\begin{description}\setlength{\itemsep}{0pt}
\item[]\mbox{ }
\item[{\rm \textsc{BaseCaseChecker}}]
\item[\emph{Input:}] $(\alpha, \beta) \in S$, and nonnegative integers $d(\alpha,\beta)$ and $l$ with $l \ge 3$
\item[\emph{Output:}] either \texttt{success} or \texttt{failure}
\item[]\mbox{ }
\item[{\rm Step 0 ($B'_0$):}] 
If $\tilde{\Delta}(l, n) > f_{\alpha,\beta}(l,n)$ holds for some pair $(l,n)$ with $n<n_L(l)$, 
then output \texttt{failure} and stop. Otherwise, go to Step~1 if $l+1 < d(\alpha,\beta)$, go to Step~2 if $d(\alpha,\beta) < 2^{2\alpha+1}$, and output \texttt{success} otherwise. 
\item[{\rm Step 1 ($B'_1$):}] 
If $\tilde{\Delta}(d, n) > f_{\alpha,\beta}(d,n)$ holds for some pair $(d, n)$ with $l+1 \le d < d(\alpha,\beta)$ and $2d \le n<n_L(d)$, 
then output \texttt{failure} and stop. Otherwise, go to Step~2. 
\item[{\rm Step 2 ($B_2$):}] 
If $\tilde{\Delta}(d, n) > f_{\alpha,\beta}(d,n)$ holds for some $(d, n)$ 
with $d(\alpha,\beta) \le d < 2^{2\alpha+1}$ and $2d \le n<d+2^{2\alpha+1}$, then output \texttt{failure} and stop. 
Otherwise, output \texttt{success}. 
\item[{\rm (End)}]
\item[]\mbox{ }
\end{description}

\begin{rem}
The parameter $l$ can be excluded from the list of inputs by adding an outer-loop for $l$; 
i.e., starting from $l=3$, if \textsc{BaseCaseChecker} outputs \texttt{success}, then we are done; otherwise, incrementing $l$ by one, we feed it to \textsc{BaseCaseChecker} and repeat the same procedure. 
\end{rem}

\begin{rem}
The computation of $d(\alpha,\beta)$ is not included in the procedure, and hence should be done in advance; 
see Claim~\ref{clm:1} for the sufficient condition for $d(\alpha,\beta)$, and also Remark~\ref{rem:dab} and Section~\ref{sec:NE} for how to compute $d(\alpha,\beta)$ in practice based on the condition. 
\end{rem}

\begin{prop}\label{prop:1}
If \textsc{BaseCaseChecker} outputs {\rm \texttt{success}}, then $\Delta(d,n) \le f_{\alpha,\beta}(d,n)$ for $n \ge d \ge l$. 
\end{prop}

\subsection{Correctness: proof of Lemma~\ref{lem:validInductive}}\label{sec:correctness}

Recall that Lemma~\ref{lem:validInductive} states that 
for given $(\alpha, \beta) \in S$, 
there exists $d(\alpha,\beta)$ such that $d \ge d(\alpha,\beta)$, $n \ge 2d$, and $n \ge d+2^{2\alpha+1}$ imply Inequality $(\ref{eqn:inductive})$. 
Recall that Inequality $(\ref{eqn:inductive})$ is 
\begin{align*}
f_{\alpha,\beta}(d-1,n-1)+2f_{\alpha,\beta}\left(d, \left\lfloor \frac{n}{2} \right\rfloor \right)+2 \le f_{\alpha,\beta}(d,n).
\end{align*}
\noindent
Since $n \ge 2d$, we have either $\left\lfloor \frac{n}{2} \right\rfloor = d$ or $\left\lfloor \frac{n}{2} \right\rfloor > d$. 

Suppose that $\left\lfloor \frac{n}{2} \right\rfloor= d$. 
In this case, the second term $2f_{\alpha,\beta}\left(d, \left\lfloor \frac{n}{2} \right\rfloor \right)$ of the left hand side of Inequality $(\ref{eqn:inductive})$ vanishes.  
As a matter of fact, a proof for the corresponding inequality immediately follows from the proof for $\left\lfloor \frac{n}{2} \right\rfloor > d$ which will be given in the following. 

Suppose that $\left\lfloor \frac{n}{2} \right\rfloor > d$. 
Observe that in general, $a^{\log_2b} = b^{\log_2a}$ for $a,b > 0$. 
Hence, letting $d(\alpha,\beta)$ be sufficiently large so that $\beta+\frac{d}{\alpha}>0$, since $n \ge d+2^{2\alpha+1}$, 
\begin{align*}
f_{\alpha,\beta}(d,n)=(n-d)^{\log_2\left(\beta+d/\alpha \right)}=\left(\beta+\frac{d}{\alpha}\right)^{\log_2\left(n-d \right)}.
\end{align*}
Using this, $(\ref{eqn:inductive})$ can be rewritten as
\begin{equation}\label{eq:target}
\left(\beta+\frac{d-1}{\alpha} \right)^{\log_2\left(n-d \right)}+2\left(\beta+\frac{d}{\alpha} \right)^{\log_2\left(\left\lfloor n/2 \right\rfloor -d \right)}+2 
\le \left( \beta+\frac{d}{\alpha} \right)^{\log_2\left(n-d \right)}.
\end{equation}
\noindent 
Note that by the integrality of $\left\lfloor \frac{n}{2} \right\rfloor$ and $d$, we have $\left\lfloor \frac{n}{2} \right\rfloor - d \ge 1$. 
It is easily seen that 
\begin{align*}
\left(\beta+ \frac{d}{\alpha} \right)^{\log_2\left(\left\lfloor  n/2 \right\rfloor -d \right)} 
&\le \left( \beta+\frac{d}{\alpha} \right)^{\log_2\left( n/2-d/2 \right)} 
= \left( \beta+\frac{d}{\alpha} \right)^{-1+\log_2(n-d)}. 
\end{align*}
Now, observe that $n \ge d + 2^{2\alpha+1}$ implies $\log_2(n-d) \ge 2\alpha+1$, 
and also that for $d \ge 2$, 
\begin{align*}
0 < \left( 1- \frac{\frac{1}{\alpha}}{\beta+\frac{d}{\alpha}} \right)=\left( 1-\frac{1}{d+\alpha\beta} \right)<1,
\end{align*}
because $(\alpha, \beta) \in S$ implies that $\alpha \ge 1$ and $\beta \ge 0$. 
Then, the left-hand side of $(\ref{eq:target})$ is bounded from above by 
\begin{align} 
		&\left( \beta+\frac{d-1}{\alpha} \right)^{\log_2(n-d)} + 2\left( \beta+\frac{d}{\alpha} \right)^{-1+\log_2(n-d)}+2  \nonumber \\
= 		&f_{\alpha, \beta}(d,n) \left[ \left( 1-\frac{\frac{1}{\alpha}}{\beta+\frac{d}{\alpha}} \right)^{\log_2(n-d)} + \frac{2}{\beta+\frac{d}{\alpha}}\right] +2 \nonumber \\
\le 	&f_{\alpha, \beta}(d,n) \left[ \left( 1-\frac{\frac{1}{\alpha}}{\beta+\frac{d}{\alpha}} \right)^{2\alpha+1} + \frac{2}{\beta+\frac{d}{\alpha}}\right] +2 \nonumber \\
\le 	&f_{\alpha, \beta}(d,n) \left[ 1- \frac{2}{f_{\alpha, \beta}(d,n)} \right] +2 \nonumber \\ 
= 		&f_{\alpha, \beta}(d,n), \nonumber
\end{align}
where the second inequality follows from Claim~\ref{clm:1} below. 

\begin{clm}\label{clm:1}
For a given $(\alpha, \beta) \in S$, there exists $d(\alpha,\beta)$ such that $d \ge d(\alpha,\beta)$, $n \ge 2d$, and $n \ge d+2^{2\alpha+1}$ imply 
\begin{equation*} 
\left( 1-\frac{\frac{1}{\alpha}}{\beta+\frac{d}{\alpha}} \right)^{2\alpha+1} + \frac{2}{\beta+\frac{d}{\alpha}} \le 1- \frac{2}{f_{\alpha, \beta}(d,n)}. 
\end{equation*}
\begin{proof}
Since $n \ge d+2^{2\alpha+1}$ implies $\log_2(n-d) \ge 2\alpha+1$, 
we have $f_{\alpha, \beta}(d,n) \ge \left( \beta+\frac{d}{\alpha} \right)^{2\alpha+1}$, and hence it suffices to show that 
\begin{equation}\label{eq:clm1a}
\left( 1-\frac{\frac{1}{\alpha}}{\beta+\frac{d}{\alpha}} \right)^{2\alpha+1} + \frac{2}{\beta+\frac{d}{\alpha}} + 2\left(\frac{1}{\beta+\frac{d}{\alpha}} \right)^{2\alpha+1} \le 1. 
\end{equation} \noindent
Letting $D=\beta+\frac{d}{\alpha}$, the left-hand side of $(\ref{eq:clm1a})$ can be rewritten as 
\[ 
\left( 1-\frac{\frac{1}{\alpha}}{D} \right)^{2\alpha+1} + \frac{2}{D} + 2\left(\frac{1}{D} \right)^{2\alpha+1} = 
1 + \frac{{\displaystyle \sum_{k=0}^{2\alpha} c(k)} D^k}{D^{2\alpha+1}}, \]
where $c(1), c(2), \ldots, c(2\alpha)$ are coefficients independent from $D$. 
In particular, the coefficient $c(2\alpha)$ of the term of maximum degree with respect to $D$ is strictly negative:   
\[ c(2\alpha) = \binom{2\alpha+1}{1} 
\cdot \left( -\frac{1}{\alpha} \right) + 2 = -2-\frac{1}{\alpha}+2 = -\frac{1}{\alpha}<0.\] 
Therefore, when $D$ is sufficiently large, the numerator $\sum_{k=0}^{2\alpha} c(k) D^k$ is strictly negative. 
Since $\alpha > 0$, one can conclude that there exists $d(\alpha,\beta)$ satisfying the desired condition, which completes the proof of Claim~\ref{clm:1}. 
\end{proof}
\end{clm}

\begin{rem}[How to calculate $d(\alpha,\beta)$ in practice]\label{rem:dab}
To compute $d(\alpha,\beta)$ satisfying the conditions of Claim~\ref{clm:1}, it is enough to determine the largest root $D^*$ of the numerator $f(D) = \sum_{k=0}^{2\alpha} c(k) D^k$ 
and simply set $d(\alpha,\beta)=\left\lceil \alpha(D^*-\beta) \right\rceil$. 
In this paper, we compute an upper bound on $\left\lceil \alpha(D^*-\beta) \right\rceil$ by elementary calculus; see Section~\ref{sec:NE} for the details. 
\end{rem}

\section{Numerical Examples}\label{sec:NE}

This section explains how \textsc{BaseCaseChecker} works using the cases $(\alpha,\beta) \in \{ (2,0), (4,0), (8,16)\}$, which yield the inequalities stated in Theorem~\ref{thm:main}. 

\subsection{Case $(\alpha,\beta)=(2,0)$}

As indicated in the proof of  Claim~\ref{clm:1}, it suffices to find $d(2,0)$ such that $d \ge d(2,0)$ implies Inequality $(\ref{eq:clm1a})$ with $(\alpha,\beta)=(2,0)$, i.e.,  
\begin{equation}\label{eq:2_0pre}
\left( 1-\frac{\frac{1}{2}}{\frac{d}{2}} \right)^{5} + \frac{2}{\frac{d}{2}} + 2\left(\frac{2}{d} \right)^{5} \le 1. 
\end{equation}
\begin{obs}\label{obs:2_0}
Inequality $(\ref{eq:2_0pre})$ holds for $d \ge 10$, hence, we can set $d(2,0):=10$. 
\begin{proof}
Observe that Inequality $(\ref{eq:2_0pre})$ is equivalent to $(d-1)^5 +4d^4+2\cdot 2^5 \le d^5$, which can be rewritten as 
\begin{equation}\label{eq:2_0}
-d^4 +10d^3-10d^2+5d+63 \le 0. 
\end{equation}
For $d \ge 10$, 
\begin{itemize}\setlength{\itemsep}{0pt}
\item
$-d^4 + 10d^3 \le -10d^3 + 10d^3 \le 0$,
\item
$-10d^2+5d \le -100d + 5d \le -95d$. 
\end{itemize}
Therefore, for $d \ge 10$, the left-hand side  of $(\ref{eq:2_0})$ is bounded from above by 
\[ -10d^3+10d^3-100d+5d+63 \le -95d+63, \]
which is negative  for $d \ge 10$. 
\end{proof}
\end{obs}

For the case $(\alpha,\beta)=(2,0)$, $l$ must be at least  four because the exponent $\log_2 \left( \frac{d}{2} \right)$ is smaller than $1$ for $d \le 3$. 
It was verified that for each of $d$ with $d \in \{4,5,6\}$, there exists a pair $(d,n)$ such that $n<n_L(d)$ while  $\tilde{\Delta}(d, n) > f_{\alpha,\beta}(d,n)$: 
\begin{align*}
&f_{2,0}(4,8) =(8-4)^{\log_2\left(4/2\right)}=4.00 < 6.00 = \tilde{\Delta}(4, 8)\\
&f_{2,0}(5,10) =(10-5)^{\log_2\left(5/2\right)} < 8.40 < 9.00 = \tilde{\Delta}(5, 10)\\ 
&f_{2,0}(6,24) =(24-6)^{\log_2\left(6/2\right)} < 97.63 < 98.00 = \tilde{\Delta}(6, 24)
\end{align*}
Hence, our approach cannot ensure $\Delta(d,n) \le f_{\alpha,\beta}(d,n)$ for $d\le 6$ although it can be true. 
On the other hand, \textsc{BaseCaseChecker} outputs {\tt success} for $l=7$.  
In what follows, we provide a few details.

\paragraph{Execution results on $(B'_0)$}
\begin{figure}[tb]
\centering
\includegraphics[height=240pt]{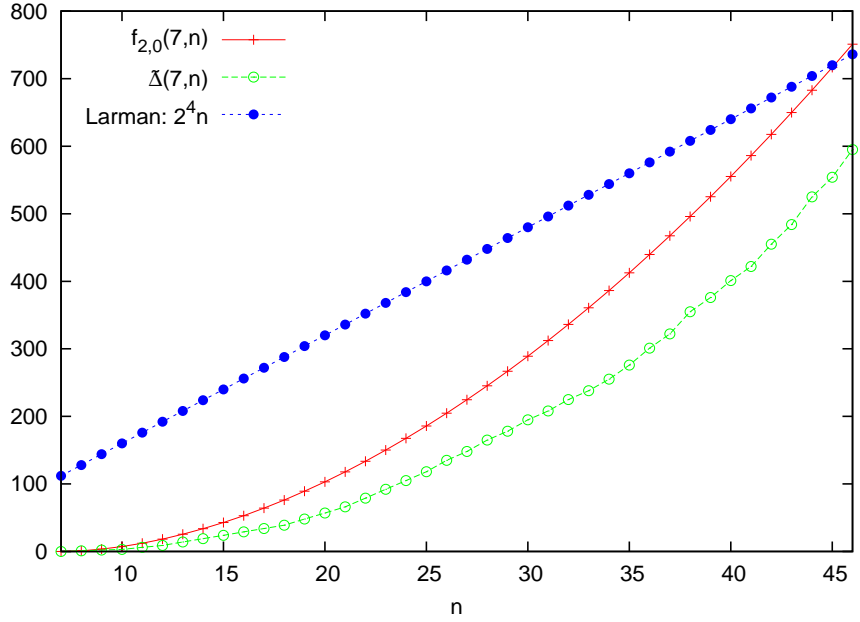}
\caption{Values of $f_{2,0}(7,n)$, $\tilde{\Delta}(7,n)$, and the generalized Larman bound for $d=7$, i.e, $2^4n$~\label{fig:12}}
\end{figure}
Figure~\ref{fig:12} shows the values of $f_{2,0}(7,n)$, $\tilde{\Delta}(7,n)$, and the generalized Larman bound for $d=7$. 
As we see from Figure~\ref{fig:12}, it was verified that for $d =7 \le n \le 45$, we have $\tilde{\Delta}(7,n) \le f_{2,0}(7,n)$, which implies $\Delta(7,n) \le f_{2,0}(7,n)$. 
Also, for $n=46$, the value of $f_{2,0}(7,n)$ is at most $750.96$ while that of the generalized Larman bound is $736$, hence, $n_L(7)=46$. 
Therefore, $\Delta(d,n) \le f_{2,0}(d,n)$ for $d = 7$, namely, $(B'_0)$ holds. 

\paragraph{Execution results on $(B'_1)$ and $(B_2)$}
Furthermore, it was verified that $(B'_1)$ and $(B_2)$ hold, where $n_L(8)=47$ and $n_L(9)=51$ when verifying $(B'_1)$, and 
$d$ ranges from $10$ to $2^{2 \cdot 2+1} = 32$ when verifying $(B_2)$. 
To sum up, by Proposition~\ref{prop:1}, 
\[ \Delta(d,n) \le f_{2,0}(d,n) = (n-d)^{\log_2\left(d/2 \right)}=(n-d)^{-1+\log_2 d} \mbox{ for } n \ge d \ge 7, \]
which yields Inequality (a) stated in Theorem~\ref{thm:main}. 

\subsection{Case $(\alpha,\beta)=(4,0)$}

In this case, it suffices to find $d(4,0)$ such that $d \ge d(4,0)$ implies
\begin{equation}\label{eq:4_0pre}
\left( 1-\frac{\frac{1}{4}}{\frac{d}{4}} \right)^{9} + \frac{2}{\frac{d}{4}} + 2\left(\frac{1}{\frac{d}{4}} \right)^{9} \le 1.
\end{equation}
\begin{obs}\label{obs:4_0}
Inequality $(\ref{eq:4_0pre})$ holds for $d \ge 36$, hence, we can set $d(4,0):=36$. 
\begin{proof}
Observe that Inequality $(\ref{eq:4_0pre})$ is equivalent to $(d-1)^9 +8d^8+2\cdot 4^9 \le d^9$, which can be rewritten as 
\begin{equation}\label{eq:4_0}
-d^8 +36d^7-84d^6+126d^5-126d^4+84d^3-36d^2+9d+524287 \le 0. 
\end{equation}
For $d \ge 36$, 
\begin{itemize}\setlength{\itemsep}{0pt}
\item
$-d^8 +36d^7 \le -36d^7 +36d^7 \le 0$,
\item
$-84d^6 +126d^5 \le -3024d^5 +126d^5 = -2898d^5$,
\item
$-126d^4 +84d^3 \le -4536d^3 +126d^3 = -4452d^3$,
\item
$-36d^2 +9d \le -1296d +9d = -1287d$. 
\end{itemize}
Therefore, for $d \ge 36$, the left-hand side  of $(\ref{eq:4_0})$ is bounded from above by 
\[ -2898d^5 -4452d^3 -1287d + 524287, \]
which is strictly negative for $d \ge 36$. 
\end{proof}
\end{obs}

\noindent
Since $l \ge d(4,0)=36$, \textsc{BaseCaseChecker} skips Step~1 and goes to Step~2 after verifying in Step~0 that $(B'_0)$ holds with $l=37$, where $n_L(37)=42946$. 
\textsc{BaseCaseChecker} verified that $(B_2)$ also holds and outputs \texttt{success}, which implies that 
\[ \Delta(d,n) \le f_{4,0}(d,n) = (n-d)^{\log_2\left( d/4 \right)}=(n-d)^{-2+\log_2 d} \mbox{ for } n \ge d \ge 37. \]
This is Inequality (b) stated in Theorem~\ref{thm:main}. 

\subsection{Case $(\alpha,\beta)=(8,16)$}

As a matter of fact, \textsc{BaseCaseChecker} runs out of computational memory for the case $(\alpha,\beta)=(8,0)$. 
This is because $n_L(d)$ exceeded the limitation on the array length in our circumstance. 

For achieving an upper bound with $\alpha=8$, one can increase the value of $\beta$. 
This makes $n_L(d)$ relatively {\em small}. 
For example, \textsc{BaseCaseChecker} outputs {\tt success} for the case $(\alpha,\beta)=(8,16)$ with $l=4$. 
It is not difficult to see that we can set $d(8,16)=8$; see Appendix~\ref{append:C} for the details. 
Since for $d \le 3$, 
\[ \Delta(d,n) \le n-d \le (n-d)^{4} \le (n-d)^{\log_2 \left( 16+ d/8 \right)} = f_{8,16}(d,n), \]
we conclude that
\[ \Delta(d,n) \le f_{8,16}(d,n) = (n-d)^{\log_2\left( 16+d/8 \right)}=(n-d)^{-3 + \log_2 d+\mO\left( 1/d \right)} \mbox{ for } n \ge d \ge 1, \]
which yields Inequality (c) stated in Theorem~\ref{thm:main}. 

\begin{rem}
A further improved upper bound of the form $f_{\alpha,\beta}(d,n)$ with $\alpha > 8$ may be proven by making $\beta$ larger. 
The resulting upper bound is, however, still in the form of $(n-d)^{\log_2 \mO(d)}$. 
\end{rem}

\begin{rem}
Since our proof method is based on only the Kalai-Kleitman inequality and the generalized Larman bound, 
one can easily apply it to a more generalized setting where we have the similar results; 
see, e.g., \cite{GaKi16} who proved an improved upper bound on the diameter of normal simplicial complexes by extending the proof of \cite{SuKi15}, a special case of this study. 
\end{rem}

\begin{rem}\label{rem:poly}
Although not surprising, our approach cannot yield any polynomial bound. 
Specifically, for an arbitrarily given polynomial function $p(d,n)$, there are infinitely many pairs $(d,n)$ such that 
Inequality $(\ref{eqn:inductive})$, the inequality which needs to shown in the inductive step, does not hold. 
\begin{proof}
See Appendix~\ref{appendix:D}. 
\end{proof}
\end{rem}

\section*{Acknowledgement}
\noindent
The author would like to thank Antoine Deza and anonymous referees for their helpful and constructive comments, which substantially improved the exposition. 
This work was supported by JSPS KAKENHI Grant Number 15H06617.

\appendix
\section{A C Code for \textsc{BaseCaseChecker} and Its Execution Results}\label{append:A}

In what follows, we show our code for \textsc{BaseCaseChecker} in the programming language C. 
The values of the parameters are those used for the case $(\alpha,\beta)=(2,0)$. Note that $d(2,0)=10$. 
It accepts the value of $l$ from the standard input. 

\begin{description}
\item[A C Code for \textsc{BaseCaseChecker}]\mbox{ }
\begin{spacing}{0.8}
{\footnotesize
\begin{verbatim}
#include<stdio.h>
#include<math.h>
#include<stdlib.h>

#define A 2 //alpha
#define B 0 //beta
#define D_AB 10 //d(alpha,beta)
#define N 1000000 //array length

double T[N]={0}, U[N]={0}; //T: tilde_D(d-1,n), U: tilde_D(d,n)
double bound_Larman(int,int); //return n*2^{d-3}
double bound_Ours(int,int); //return f_{A,B}(d,n)

void initialize(void); //initialize T and U
void Update(int); //update U using T
void check(int,int); //check if tilde_D(d,n)<=f(d,n)

int main(void)
{
   int d,n,d_max=(int)pow(2.0,2.0*A+1); //d, n, and the maximum of d
   int flag = 0; //takes 1 if tilde_D(d,n)>f(d,n) holds
   initialize();
   
   int l=3;
   printf("Enter l (>=3): ");
   scanf("%d",&l);
   
   //Compute tilde_D(d,n) for d=l
   for(d=3; d<l; d++){
      Update(d+1);
   }
   
   //B0
   n=l; printf("\n");
   while(n<N && bound_Larman(d,n)>bound_Ours(d,n)){
      check(d,n);
      n++;
   }
   printf("- n_L(%d) = %d\n",d,n);
   printf("(B0) OK\n");
    
   Update(d+1);
   d++;
   
   //B1
   while(d<D_AB){
      n=2*d;
      while(bound_Larman(d,n)>bound_Ours(d,n)){
         check(d,n);
         n++;
      }
      printf("- n_L(%d) = %d\n",d,n);
      Update(d+1);
      d++;
   }
   printf("(B1) OK\n");
   
   //B2
   while(d<d_max){
      n=2*d;
      int count = 0;
      while(n<d+d_max){
         check(d,n);
         n++;
         count++;
      }
      printf("- # pairs (%d,n) checked = %d\n",d,count);
      Update(d+1);
      d++;
   }
   printf("(B2) OK\n");
   printf("\n****** SUCCESS ******");
   
   return 0;
}

double bound_Larman(int d, int n){return n*pow(2.0,d-3);}
double bound_Ours(int d, int n){return pow(1.0*(n-d),log(1.0*d/A+B)/log(2));}

void initialize(void)
{
   int i;
   for(i=0; i<N; i++){
      int n = i+3; //i = n - 3
      U[i] = 1.0*(n-3); //use the Hirsch bound for d=3
   }
}

void Update(int d)
{
   int i;
   for(i=0; i<N; i++) T[i] = U[i];
   for(i=0; i<N; i++){
      int n = i+d; //i = n - d
      if(n < 2*d){U[i] = T[i];} //tilde_D(d,n)=tilde_D(d-1,n-1)
      else{U[i] = T[i]+2*U[n/2-d]+2;} //tilde_D(d,n)=tilde_D(d-1,n-1)+2*tilde_D(d,n/2)+2
   }
}

void check(int d,int n)
{
   int i=n-d;
   if(bound_Ours(d,n)<U[i]){
      printf("Error: %.1f [Ours] < %.1f [tilde] (%d,%d)\n",bound_Ours(d,n),U[i],d,n); 
      printf("\n****** FAILURE ******"); exit(1);
   }else if(n==N-1){
      printf("Error: Out of Memory\n"); 
      printf("\n****** FAILURE ******"); exit(1);
   }
}
\end{verbatim}
}
\end{spacing}
\end{description}

\noindent
Next, we show some execution results of the above code. 
\begin{description}
\item[Case $(\alpha,\beta)= (2,0)$ with $l=7$ and $d(2,0)=10$]\mbox{ }
\begin{spacing}{0.8}
{\footnotesize
\begin{verbatim}
Enter l (>=3): 7

- n_L(7) = 46
(B0) OK
- n_L(8) = 47
- n_L(9) = 51
(B1) OK
- # pairs (10,n) checked = 22
- # pairs (11,n) checked = 21
- # pairs (12,n) checked = 20
//snip
- # pairs (29,n) checked = 3
- # pairs (30,n) checked = 2
- # pairs (31,n) checked = 1
(B2) OK

****** SUCCESS ******
\end{verbatim}
}
\end{spacing}

\item[Case $(\alpha,\beta)= (2,0)$ with $l=6$ and $d(2,0)=10$]\mbox{ }
\begin{spacing}{0.8}
{\footnotesize
\begin{verbatim}
Enter l (>=3): 6

Error: 97.6 [Ours] < 98.0 [tilde] (6,24)

****** FAILURE ******
\end{verbatim}
}
\end{spacing}

\item[Case $(\alpha,\beta)= (4,0)$ with $l=37$ and $d(4,0)=36$]\mbox{ }
\begin{spacing}{0.8}
{\footnotesize
\begin{verbatim}
Enter l (>=3): 37

- n_L(37) = 42946
(B0) OK
(B1) OK
- # pairs (38,n) checked = 474
- # pairs (39,n) checked = 473
- # pairs (40,n) checked = 472
//snip
- # pairs (509,n) checked = 3
- # pairs (510,n) checked = 2
- # pairs (511,n) checked = 1
(B2) OK

****** SUCCESS ******
\end{verbatim}
}
\end{spacing}

\item[Case $(\alpha,\beta)= (4,0)$ with $l=36$ and $d(4,0)=36$]
\begin{spacing}{0.8}
{\footnotesize
\begin{verbatim}
Enter l (>=3): 36

Error: 1469828390203.3 [Ours] < 1469922992914.0 [tilde] (36,6928)

****** FAILURE ******
\end{verbatim}
}
\end{spacing}

\end{description}

\section{The computation of $d(8,16)$}\label{append:C}

It suffices to find $d(8,16)$ such that $d \ge d(8,16)$ implies 
\begin{equation}\label{eq:last_pre}
\left( 1-\frac{\frac{1}{8}}{16+\frac{d}{8}} \right)^{17} + \frac{2}{16+\frac{d}{8}} + 2\left(\frac{1}{16+\frac{d}{8}} \right)^{17} \le 1. 
\end{equation}
For notational simplicity, set $D:=\frac{d}{8}+16$, and rewrite Inequality (\ref{eq:last_pre}) as 
\begin{equation}\label{eq:last}
\left( D-\frac{1}{8} \right)^{17} + 2D^{16} + 2 \le D^{17}.
\end{equation}
We prove that Inequality (\ref{eq:last}) holds for $D \ge 17$. 
If this is true, then Inequality (\ref{eq:last_pre}) is satisfied for $\frac{d}{8}+16 \ge 17$, i.e., for $d \ge 8$. 
By simple calculus, Inequality (\ref{eq:last}) is rewritten as 
\begin{align*}
&D^{17} -\frac{17}{8}D^{16} +2D^{16} +2 -D^{17} + \sum_{i=2}^{17} \binom{17}{i} D^{17-i} \left( -\frac{1}{8} \right)^i \le 0\\
\iff &-\frac{1}{8}D^{16}+\binom{17}{2}\frac{1}{8^2}D^{15}-\binom{17}{3}\frac{1}{8^3}D^{14}+\binom{17}{4}\frac{1}{8^4}D^{13}- 
\cdots -\binom{17}{16}\frac{1}{8^{16}}D-\binom{17}{17}\frac{1}{8^{17}}+2\le 0. 
\end{align*}
We observe that for $D \ge 17$, 
\begin{align}
-\frac{1}{8}D^{16}+\binom{17}{2}\frac{1}{8^2}D^{15} \ \le& \ -\frac{17}{8}D^{15}+\frac{17\cdot 16}{2\cdot 1}\frac{1}{8^2}D^{15} = 0,\nonumber \\ 
-\binom{17}{3}\frac{1}{8^3}D^{14}+\binom{17}{4}\frac{1}{8^4}D^{13}  \ \le& \ D^{13} \frac{1}{8^3} \binom{17}{3} \left( -17 + \frac{14}{4}\cdot \frac{1}{8}\right) \le 0, \label{tiny}\\
-\binom{17}{5}\frac{1}{8^5}D^{12}+\binom{17}{6}\frac{1}{8^6}D^{11}  \ \le& \ D^{11} \frac{1}{8^5} \binom{17}{5} \left( -17 + \frac{12}{6}\cdot \frac{1}{8}\right)\le 0,\nonumber \\
 \vdots & \nonumber \\
-\binom{17}{15}\frac{1}{8^{15}}D^2 + \binom{17}{16}\frac{1}{8^{16}}D  \ \le& \ D \frac{1}{8^{15}} \binom{17}{15} \left( -17 + \frac{2}{16}\cdot \frac{1}{8}\right)\le 0. \nonumber 
\end{align}
We remark that here, the tiny term of $-\binom{17}{17}\frac{1}{8^{17}}+2$, which is at most $2$, was ignored. 
It can be included, for example, in Inequality {\rm (\ref{tiny})} because for $D \ge 17$, 
\[ D^{13} \frac{1}{8^3} \binom{17}{3} \left( -17 + \frac{14}{4}\cdot \frac{1}{8}\right) \le 
17^{13} \frac{1}{17^3} \left( -1 \right)  \ll -2. \]
This completes the proof. 

\section{Proof of Remark~\ref{rem:poly}}\label{appendix:D}
We want to prove that for an arbitrarily given polynomial function $p(d,n)$, there are infinitely many pairs $(d,n)$ such that 
the inequality, which needs to be proved in the inductive step of our proof method, does not hold; i.e., 
\[p(d-1,n-1)+2p\left(d, \left\lfloor \frac{n}{2} \right\rfloor \right)+2 > p(d,n), \]
holds, even when $d$ and $n$ are sufficiently large. 

The polynomial function $p(d,n)$ can be rewritten as 
\[ p(d,n) = \sum_{i=0}^k g_i(d) n^i \]
for some nonnegative integer $k$, where $g_i(d)$ is a polynomial function of $d$ for each $i \in \{0,1, \ldots, k \}$. 

We can assume that for any $d$, $g_k(d) \ge 0$; 
otherwise $p(d,n)$ cannot be a valid upper bound on $\Delta(d,n)$ because, in this case, $p(d,n)$ is negative for sufficiently large $n$. 
In what follows, we consider the case when $d$ is larger than the maximum root of $g_k$. 
Note that in this case, $g_k(d)>0$. 
It is easily seen that
\[ \frac{p(d,n)}{g_k(d) n^k} = \sum_{i=0}^k \frac{g_i(d)}{g_k(d)} \cdot \frac{1}{n^{k-i}}=1+\sum_{i=0}^{k-1} \frac{g_i(d)}{g_kd)} \cdot \frac{1}{n^{k-i}}. \]
Therefore, for any $\epsilon>0$, when $n$ is sufficiently large, 
\[ (1-\epsilon) g_k(d) n^k \le p(d,n) \le (1+\epsilon) g_k(d) n^k. \]
By similar arguments, 
\begin{align*}
(1-\epsilon) g_k(d-1) (n-1)^k &\le p(d-1,n-1) \le (1+\epsilon) g_k(d-1) (n-1)^k,\\
(1-\epsilon) g_k(d) \left(\left\lfloor \frac{n}{2} \right\rfloor \right)^k &\le p\left(d, \left\lfloor \frac{n}{2} \right\rfloor \right) \le (1+\epsilon) g_k(d) \left(\left\lfloor \frac{n}{2} \right\rfloor \right)^k.
\end{align*} 
Using these relations, 
\begin{align*}
\frac{p(d-1,n-1)}{p(d,n)} &\ge \frac{1-\epsilon}{1+\epsilon} \cdot \frac{g_k(d-1)}{g_k(d)} \cdot  \left(1-\frac{1}{n} \right)^k
\end{align*}
Since $g_k$ is a polynomial function of $d$, the ratio of $\frac{g_k(d-1)}{g_k(d)}$ gets arbitrarily close to $1$ by taking $d$ sufficiently large. 
Therefore, the right-hand side gets arbitrarily close to $1$ by taking $d$ and $n$ sufficiently large, and $\epsilon$ sufficiently small. 

On the other hand, assuming that $n$ is even for convenience, 
\begin{align*}
\frac{p\left(d, \left\lfloor \frac{n}{2} \right\rfloor \right)}{p(d,n)} 
&\ge \frac{1-\epsilon}{1+\epsilon} \left(\frac{1}{2} \right)^k 
\end{align*}
The right-hand side is bounded from below by a positive constant factor, say, $\frac{1}{3} \left(\frac{1}{2} \right)^k$ when $\epsilon \le \frac{1}{2}$. 
To sum up, if both $d$ and $n$ are sufficiently large, then 
\[\frac{p(d-1,n-1)}{p(d,n)} + \frac{p\left(d, \left\lfloor n/2 \right\rfloor \right)}{p(d,n)} > 1, \] 
which implies $ p(d-1,n-1) + p\left(d, \left\lfloor n/2 \right\rfloor \right) +2 > p(d,n)$. 

\end{document}